\documentclass[11pt,reqno]{amsart}
\usepackage[margin=1in]{geometry}
\numberwithin{equation}{section}
\usepackage{amsmath,amssymb}
\usepackage{titlesec}
\usepackage{bm}
\usepackage{enumerate}

\titleformat{\section}[block]
{\centering\bfseries\Large}                    
{\bfseries~\thesection.}           
{1em}                                          
{}                                             

\newtheorem{thm}{Theorem}[section]

\newtheorem{lem}[thm]{Lemma}
\newtheorem{prop}[thm]{Proposition}
\theoremstyle{definition}
\newtheorem{defn}[thm]{Definition}
\newtheorem{rem}[thm]{Remark}

\newcommand{\p}{\partial}

\newcommand{\R}{\mathbb{R}}

\newcommand{\vp}{\varphi}

\newcommand{\Om}{\Omega}

\newcommand{\DN}{\Delta_\infty^N}

\title[Harnack Inequality]{Harnack Inequality for Nonlinear Equations Driven by the Normalized Infinity-Laplacian}

\author{Ahmed Mohammed and Carson Pocock}

\address{Department of Mathematical Sciences, Ball State University, Muncie, IN 47306, USA}
\email{amohammed@bsu.edu}

\address{Department of Mathematical Sciences, Ball State University, Muncie, IN 47306, USA}
\email{carson.pocock@bsu.edu}

\begin{document}

\begin{abstract}
We establish interior Harnack inequalities for nonnegative viscosity
solutions of
\[
\Delta_\infty^N u=f(u)+g(u)|Du|^q,\qquad 0<q<2,
\]
under natural monotonicity, sign, and Keller--Osserman-type assumptions on
the nonlinearities. We first prove a local Harnack inequality for positive
supersolutions of a model equation with bounded zeroth-order and gradient
coefficients; a family of exponential supersolutions shows that the range
\(0<q<2\) is sharp. We then construct radial blow-up barriers and obtain a
global estimate in terms of the distance to the boundary and the inverse of
an associated Keller--Osserman integral. Combining this estimate with a
comparison principle, positive regularization, and a Harnack-chain argument
yields the desired inequality on every connected subdomain compactly
contained in the ambient domain.
\end{abstract}

\maketitle

\section{Introduction}

The infinity-Laplacian
\[
\Delta_\infty u=\langle D^2u\,Du,Du\rangle
\]
arises formally as the limiting operator of the \(p\)-Laplacian
\[
\Delta_p u=\operatorname{div}\!\left(|Du|^{p-2}Du\right)
\]
as \(p\to\infty\); see \cite{BDM}. It was introduced by Aronsson
\cite{ARO} in connection with variational problems in \(L^\infty\) and
absolutely minimizing Lipschitz extensions. Such extensions formally satisfy
\(\Delta_\infty u=0\); comprehensive accounts of this theory can be found in
\cite{ACJ, MCR}. The viscosity-solution framework of Crandall and Lions
\cite{CRL} provides the appropriate interpretation of this highly
degenerate equation and supplies the comparison and stability tools used in
its analysis.

A major conceptual advance was the tug-of-war interpretation of the
normalized infinity-Laplacian developed by Peres, Schramm, Sheffield, and
Wilson \cite{PSSW}. For the Dirichlet problem
\begin{equation}\label{eq:homogeneous}
	\Delta_\infty^N u = 0 \quad \text{in } \Omega, \qquad
	u = b \quad \text{on } \partial\Omega,
\end{equation}
where
\[
\Delta_\infty^N u:=\frac{\Delta_\infty u}{|Du|^2},
\]
the game values satisfy a dynamic programming principle and converge
uniformly to the viscosity solution. Complementing this probabilistic
viewpoint, Lu and Wang \cite{LUW-2} developed a PDE approach to the
nonhomogeneous equation
\[
\Delta_\infty^N u = f(x),
\]
establishing comparison, existence, and uniqueness results. Related
Dirichlet problems for equations of the form
\(\Delta_\infty u=f(x,u)\) were studied in \cite{TA-1, TA-2}.

Harnack inequalities provide quantitative control of the local oscillation
of positive solutions and are a central tool in elliptic theory. For
\(\infty\)-harmonic functions, the first Harnack inequality was obtained by
Manfredi and Lindqvist \cite{MAL} by passing to the limit in the
\(p\)-harmonic inequality. Bhattacharya \cite{TBH} later gave a direct proof
for nonnegative infinity-superharmonic functions. This barrier method was
adapted in \cite{BIV} to infinity-Laplacian equations with lower-order terms
depending on the solution and its gradient.

\vspace{.2cm}

In this paper we consider nonnegative viscosity solutions of
\begin{equation}\label{me}
	\Delta_\infty^N u = f(u) + g(u)\,|D u|^{q}
	\quad \text{in } \Omega,
\end{equation}
where \(0<q<2\) and \(f,g\) are continuous nondecreasing functions satisfying
the structural assumptions stated below. The main results may be summarized
as follows.

\begin{enumerate}
\item We first establish a local Harnack inequality for positive viscosity
supersolutions of
\[
\Delta_\infty^N u=A(x)u+B(x)|Du|^q|u|^{1-q},
\]
when \(A\) and \(B\) are bounded above. The estimate is uniform on
sufficiently small balls. Exponential supersolutions show that no analogous
universal estimate can hold for \(q\ge2\), so the exponent range in the paper
is sharp; see Remark~\ref{sharp-q}.

\item Under the Keller--Osserman condition introduced below, radial solutions
of an associated initial-value problem provide blow-up barriers. They yield a
global a priori estimate
\[
u(x)\le \mathcal Q(d_\Omega(x)),
\]
where \(\mathcal Q\) is expressed through the inverse \(\Phi\) of a
Keller--Osserman integral. Under the quasi-monotonicity condition
{\bfseries(Q-m)}, this improves to
\[
u(x)\le\Phi(d_\Omega(x)).
\]

\item Combining the global estimate with a comparison principle, the
regularization \(u+\varepsilon\), and a finite Harnack chain, we prove an
interior Harnack inequality on every connected open set compactly contained
in \(\Omega\). For \(0<q<1\), the hypotheses are formulated through
{\bfseries(Q-M)}. For \(1\le q<2\), they use {\bfseries(Q-m)} and
{\bfseries(KO)$_q$}, together with {\bfseries(Q-q)} when \(q>1\).
\end{enumerate}

The proof therefore combines local barriers of infinity-harmonic type with
global Keller--Osserman bounds tailored to the nonlinear absorption and
gradient terms in \eqref{me}.

\vspace{.2cm}

Below we introduce notational conventions that will be used throughout the paper.

\begin{itemize}
	\item o stands for the origin in $\mathbb R^n$.
	
	\vspace{.2cm}
	\item $\mathbb R^+_0:=[0,\infty),\;\;\;\;\mathbb R^+:=(0,\infty)$
		\vspace{.2cm}
		
		\item $B(x,r)$ is the open ball in $\mathbb R^n$ of radius $r>0$ and centered at $x$.

	\vspace{.2cm}
	\item $\Omega\subset\mathbb R^n$ stands for an open subset with non-empty boundary $\partial \Omega$.
	
	\vspace{.2cm}
	\item For a non-empty subset $E\subset\Omega$, we write $\textup{dist}(E,\partial\Omega):=\inf\{|x-y|:x\in E,\,y\in\partial\Omega\}.$
	
		\vspace{.2cm}
		
		\item For $x\in\Omega$, we write $d_\Omega(x):=\textup{dist}(\{x\},\partial\Omega)$.

	\vspace{.2cm}
	\item $\textup{USC}(\Omega)$ denotes the class of upper-semicontinuous functions in $\Omega$.

	\vspace{.2cm}
	\item $\textup{LSC}(\Omega)$ denotes the class of lower-semicontinuous functions in $\Omega$.
	
	\vspace{.2cm}
	\item $\mathcal C(\Omega):=\textup{USC}(\Omega)\cap \textup{LSC}(\Omega)$
	
	\vspace{.2cm}

	\item $C^2(\Omega)$ denotes the class of twice continuously differentiable functions $u:\Omega\to \mathbb R$.
	
	\vspace{.2cm}
	\item $\mathcal S^{n\times n}(\R)$ denotes the set of all $n\times n$ symmetric matrices with real entries.
\end{itemize}

\section{Main Results}

\noindent We begin by considering the equation
\begin{equation}\label{e-1}
	\DN u=A(x)u+B(x)|Du|^q|u|^{1-q},
\end{equation}
where $A, B\in \mathcal C(\Omega)$ have uniform upper bounds. In other words, we assume that there are constants $A_0>0$ and $B_0>0$ such that $A(x)\le A_0$ and $ B(x)\le B_0$ for all $x\in \Omega$.

The following result is the first step toward establishing a Harnack inequality for solutions of (\ref{me}). It establishes a Harnack inequality for positive viscosity supersolutions of (\ref{e-1}).

\begin{thm}\label{Hi-1} Let $\Om\subseteq\R^n$ be an open set, and let $0<q<2$. Assume that
	\[
	A(x)\le A_0,\qquad B(x)\le B_0
	\qquad\text{for all }x\in\Omega,\]
	for some positive constants $A_0,B_0$. Given $\alpha>1$, there is a constant $r_0>0$, depending only on $\alpha, q$, $|A_0|$ and $|B_0|$, such that for any ball $B(x_0,2r)\subset\Om$ with $0<r<r_0$ and every positive viscosity supersolution $u\in\textup{LSC}(\Om)$ of (\ref{e-1}), we have
	\begin{equation}\label{har-1}
		\sup_{B(x_0,r/3)}u\le 3^\alpha\inf_{B(x_0,r/3)}u
	\end{equation}
\end{thm}

\vspace{.3cm}
The exponent range $0<q<2$ is sharp; see
Remark~\ref{sharp-q} in Appendix~\ref{app:sharpness}.

In order to state our other main result, we need to discuss some relevant growth conditions on the continuous functions $f,g:\mathbb R\to\mathbb R$ that appear in equation (\ref{me}). We begin with a condition needed to establish a suitable comparison principle for (\ref{me}).

\begin{description}
	\item[$\bm{(\mathcal{P})}$] $f,g:\R\to\R$ are nondecreasing continuous functions such that 
	\begin{enumerate}[\bfseries (a)]
		
		\vspace{.1cm}
		\item  $f(-t)<0<f(t)$ for all $t>0$,

		\vspace{.2cm}
		
		\item either $f$ or $g$ is strictly increasing.

	\end{enumerate}
\end{description}
\vspace{.2cm}

Note that condition $\bm{(\mathcal P)\,(\mathbf a)}$ implies that $f(0)=0$.
Whenever, in addition to $\bm{(\mathcal P)}$, we assume that
\[
g(t)\ge 0 \qquad \text{for all } t>0,
\]
we shall denote this strengthened condition by $\bm{(\mathcal P^+)}$.

A uniform global upper bound on all viscosity subsolutions of (\ref{me}) will be one of the tools we will use to establish our Harnack inequality. For this we need suitable growth conditions on $f$ and $g$ at infinity. This growth condition is captured by the following integral condition, reminiscent of the classical Keller--Osserman condition:

\vspace{.2cm}
\begin{description}
	\item[(KO)$_q$] $\displaystyle{\int_1^\infty\frac{ds}{\sqrt{F(s)}+(G(s))^{\frac{1}{2-q}}}<\infty}$.
\end{description}

\vspace{.2cm}
Here, $F$ and $G$ stand for the antiderivatives of $f$ and $g$, respectively, that vanish at zero. That is
\[F(t):=\int_0^t f(s)\,ds,\;\;\text{and}\;\;\;G(t):=\int_0^t g(s)\,ds.\]

\noindent Finally, given a function $h:\mathbb R^+\to\mathbb R^+$  we consider the following conditions.
\begin{description}
	\item[{\bfseries (Q-M)}] There are constants $0<\varrho<1$ and $\nu>0$ such that 
	\[\frac{h(s)}{s^\nu}\ge \varrho\,\frac{h(t)}{t^\nu},\;\;\;s\ge t>0.\]
\end{description}

\noindent It is obvious that the condition {\bfseries(Q-M)} above implies the following condition on $h$:

\begin{description}
	\item[{\bfseries (Q-m)}]  There is a constant $0< \mu<1$ such that 
	\[h(s)\ge \mu\, h(t),\;\;\;s\ge t>0.\]
\end{description}
In fact, if {\bfseries(Q-M)} holds then {\bfseries(Q-m)} holds with $\mu:=\varrho$.

\vspace{.2cm}
\noindent For $0<q<2$, conditions {\bfseries(Q-M)} and {\bfseries(Q-m)} will be used for the following function whenever $f$ and $g$ satisfy $\bm{(\mathcal P^+)}$.
\begin{equation}\label{hfg}
	h(s):=\left(\frac{f(s)}{s}\right)^{1/2}+\left(\frac{g(s)}{s^{1-q}}\right)^{1/(2-q)},\;\;\;s>0.\tag{$\bm f$-$\bm g$}
\end{equation}

\noindent As shown in Appendix~\ref{app:growth}, condition {\bfseries(Q-M)} for the function $h$ in (\ref{hfg}) implies condition {\bfseries(KO)$_q$}.

The final condition, used when $1<q<2$, is the following.

\begin{description}
	\item[{\bfseries (Q-q)}]  One of the following limits holds.
	\[\lim_{t\to\infty}f(t)=\infty\;\;\;\;\;\;\text{or}\;\;\;\;\;\;\lim_{t\to\infty}g(t)=\infty.\]
\end{description}

\vspace{.2cm}
We are now in a position to state the Harnack inequality for nonnegative solutions of (\ref{me}) as follows.
\begin{thm}\label{mhi} Suppose $f,g:\mathbb R\to\mathbb R$ are nondecreasing and continuous functions that satisfy condition $\bm{(\mathcal P^+)}$. Let $\mathcal O$ be a connected and open set that is compactly contained in $\Omega$. 
	
\begin{enumerate}[\bfseries(1)]
	\item If $\,0<q<1$ and the function $h$ in {\bfseries(f-g)} satisfies
	condition {\bfseries(Q-M)}, then there is a positive constant $C$ that
	depends only on $q,\nu,\varrho,\mathcal O$, and
	$\operatorname{dist}(\mathcal O,\partial\Omega)$ such that for any
	nonnegative viscosity solution $u\in\mathcal C(\Omega)$ of (\ref{me}) we
	have
	\begin{equation}\label{Hi-2}\sup_{\mathcal O}u\le C\inf_{\mathcal O}u.
	\end{equation}
	
\vspace{.2cm}
\item If $1\le q<2$, the function $h$ in {\bfseries(f-g)} satisfies {\bfseries(Q-m)},
the functions $f$ and $g$ satisfy {\bfseries(KO)$_q$}, and, when $1<q<2$,
they also satisfy {\bfseries(Q-q)}, then there is a positive constant $C$,
depending on $q,\mu,f,g,\mathcal O$, and
$\operatorname{dist}(\mathcal O,\partial\Omega)$ only, such that for any
nonnegative viscosity solution $u\in\mathcal C(\Omega)$ of (\ref{me}) the
inequality (\ref{Hi-2}) holds.
\end{enumerate}

\end{thm}

\section{Preliminaries}
Because of the singular and degenerate elliptic nature of the normalized infinity-Laplacian, the appropriate framework for studying solutions of such equations is that of viscosity solutions. To recall the definition, we begin with the following notations: Given $\phi\in C^2(\Omega)$ and $x\in\Omega$, we write
\begin{align*}
	\Delta_\infty^{N,+}\phi(x)&:=\left\{\begin{array}{ll}
		|D\phi(x)|^{-2}\left\langle D^2\phi(x)D\phi(x)\,,\,D\phi(x)\right\rangle&\text{if}\;\;D\phi(x)\not=o\\[.3cm]
		\max\left\{\left\langle D^2\phi(x)e,e\right\rangle:|e|=1\right\}
		&\text{if}\;\;D\phi(x)=o,\end{array}	
	\right.\\[.3cm]
	\Delta_\infty^{N,-}\phi(x)&:=\left\{\begin{array}{ll}
		|D\phi(x)|^{-2}\left\langle D^2\phi(x)D\phi(x)\,,\,D\phi(x)\right\rangle&\text{if}\;\;D\phi(x)\not=o\\[.3cm]
		\min\left\{\left\langle D^2\phi(x)e,e\right\rangle:|e|=1\right\}
		&\text{if}\;\;D\phi(x)=o.\end{array}	
	\right.
\end{align*}

When $D\phi(x)\not=o$, it is convenient to write $\Delta_\infty^N\phi(x)$ for $\Delta_\infty^{N,+}\phi(x)=\Delta_\infty^{N,-}\phi(x)$. With these notations on hand we now recall the concepts of viscosity subsolution, supersolution and solution to 
\begin{equation}\label{gee}\Delta_\infty^N u=H(x,u,Du),\;\;\;x\in\Omega,
\end{equation} 
where $H:\Omega\times \mathbb R \times\mathbb R^n\to\mathbb R$ is a continuous function. 

\begin{defn}
	
\begin{enumerate}[\bfseries(a)]
\item A function $u\in \textup{USC}(\Omega)$ is said to be a viscosity subsolution of (\ref{gee}) if, whenever $x_0\in\Omega$ and $\phi$ is $C^2$ in a neighborhood of $x_0$ such that $u-\phi$ has a local maximum at $x_0$, then
\[\Delta_\infty^{N,+}\phi(x_0)\ge H(x_0,u(x_0),D\phi(x_0)).\] 
\item A function $u\in \textup{LSC}(\Omega)$ is said to be a viscosity supersolution of (\ref{gee}) if, whenever $x_0\in\Omega$ and $\phi$ is $C^2$ in a neighborhood of $x_0$ such that $u-\phi$ has a local minimum at $x_0$, then
\[\Delta_\infty^{N,-}\phi(x_0)\le H(x_0,u(x_0),D\phi(x_0)).\] 

\item A function $u\in \mathcal C(\Omega)$ is said to be a viscosity solution of (\ref{gee}) provided that $u$ is both a viscosity subsolution and a supersolution of (\ref{gee}) in $\Omega$.
\end{enumerate}

\end{defn}

\noindent Now we turn to our main objective of studying the Harnack inequality
for solutions to equation (\ref{me}). A comparison principle is a critical
tool for this.

\vspace{.2cm}
The following comparison principle holds. We take $\varpi:\mathbb R^+_0\to\mathbb R^+_0$ to be a continuous function such that $\varpi(t)>0$ for $t>0$.

\begin{prop}[Comparison Principle]\label{cp} Let $\Omega\subset\R^n$ be a bounded open set, and suppose $f$ and $g$ satisfy condition $\bm{(\mathcal{P}^+)}$. Let $u\in \textup{USC}(\overline{\Om})$ and $v\in \textup{LSC}(\overline{\Om})$ be such that the following hold in $\Om$ in the viscosity sense:
	\begin{equation}\label{susu}
		\Delta^N_\infty u\ge f(u)+g(u)\varpi(|Du|),\;\;\;\text{and}\;\;\;\;\Delta^N_\infty  v\le f(v)+g(v)\varpi(|Dv|).
	\end{equation}
	If $u\le v$ on $\p\Om$, then $u\le v$ in $\Om$.
\end{prop}

\noindent{\bfseries Proof:} Suppose $u>v$ somewhere in $\Om$ so that $(u-v)(x_0)=\max_{\overline{\Om}}(u-v)>0$ for some $x_0\in\Omega$. Let
\[\psi_j(x,y):=u(x)-v(y)-\frac j4|x-y|^4,\;\;\;\;\;j=1,2,\ldots,\;\;\;(x,y)\in \overline{\Om}\times\overline{\Om},\] and let $(x_j,y_j)\in \overline{\Om}\times \overline{\Om}$ such that
\[\psi_j(x_j,y_j)=\max_{\overline{\Om}\times\overline{\Om}}\psi_j(x,y).\] 
Thus for $j\in\mathbb N$ we have
\begin{equation}\label{lil}u(x)-v(y)-\frac j4|x-y|^4\le u(x_j)-v(y_j)-\frac j4|x_j-y_j|^4,\;\;\;\;\;\;\forall\,(x,y)\in\overline{\Om}\times\overline{\Om}.
\end{equation}

Passing to a subsequence, if necessary, we suppose that $(x_j,y_j)\to (\overline{x},\overline{y})\in \overline{\Omega}\times\overline{\Omega}$. It is well-known that
\begin{equation}\label{boo}
	\lim_{j\to\infty}\frac j4|x_j-y_j|^4=0.
\end{equation} As a consequence we get $\overline{x}=\overline{y}$. But then, since $u,\,-v\in \textup{USC}(\overline{\Omega})$, and using (\ref{boo}) we have
\[u(\overline{x})-v(\overline{x})\ge \limsup_{j\to\infty}(u(x_j)-v(y_j))=\limsup_{j\to\infty}\psi_j(x_j,y_j)\ge u(x_0)-v(x_0)>0.\] Since $u\le v$ on $\partial\Omega$, the above inequality implies that  $(\overline{x},\overline{x})\in\Omega\times\Omega$. Consequently $(x_j,y_j)\in \Om\times\Om$ for all sufficiently large indices $j$.

Let us show that $x_j\not= y_j$ for all sufficiently large $j$. For this, first observe that
\begin{equation}\label{ll}
	u(x_j)-v(y_j)\ge u(x_j)-v(y_j)-\frac j4|x_j-y_j|^4\ge u(x_0)-v(x_0)>0.
\end{equation}

\noindent Using $y=y_j$ in (\ref{lil}) we see that
\[u(x)\le \phi_j(x):=u(x_j)+\frac j4|x-y_j|^4-\frac j4|x_j-y_j|^4\;\;\;\;\forall\,x\in \Om.\] Since $u(x_j)=\phi_j(x_j)$ and $u$ is a subsolution we see that
\begin{equation}\label{ss}
	\Delta_\infty^{N,+}\phi_j(x_j)\ge f(u(x_j))+g(u(x_j))\varpi(|D\phi_j(x_j)|).
\end{equation}

\noindent Similarly, using $x=x_j$ in (\ref{lil}) we see that
\[v(y)\ge \vp_j(y):=v(y_j)-\frac j4|x_j-y|^4+\frac j4|x_j-y_j|^4,\;\;\forall\,y\in\Omega,\;\;\text{and}\;\;\;v(y_j)=\vp_j(y_j).\]

\noindent Recalling that $v$ is a supersolution we have
\begin{equation}\label{sst}
	\Delta_\infty^{N,-}\vp_j(y_j)\le f(v(y_j))+g(v(y_j))\varpi(|D\vp_j(y_j)|).
\end{equation}

\noindent Suppose now that $x_j=y_j$.

Then we see that 
\[D\phi_j(x_j)=o=D\vp_j(y_j),\;\;\;\;\text{and}\;\;\; D^2\phi_j(x_j)=0=D^2\vp_j(y_j).\] Consequently we have  \begin{equation}\label{ze}
\Delta_\infty^{N,+}\phi_j(x_j)=0=\Delta_\infty^{N,-}\vp_j(y_j).
\end{equation} Let us first suppose that $\varpi(0)=0$. Then (\ref{ze}), together with (\ref{ss}) and (\ref{sst}),  imply that
\begin{equation}\label{opq}
	f(u(x_j))\le 0,\;\;\;\;\text{and}\;\;\;\;f(v(y_j))\ge 0.
\end{equation}

\noindent In view of $\bm{(\mathcal P)\,(\mathbf a)}$, the conclusion in
(\ref{opq}) shows that $u(x_j)\le0$ and $v(y_j)\ge0$. Therefore
$u(x_j)-v(y_j)\le0$, which contradicts (\ref{ll}). If, on the other hand,
$\varpi(0)>0$, we obtain the following from (\ref{ss}), (\ref{sst}), and
(\ref{ze}).
\[f(u(x_j))-f(v(y_j))+\varpi(0)(g(u(x_j))-g(v(y_j)))\le 0.\]
However, this also contradicts $\bm{(\mathcal P)\,(\mathbf b)}$ and
(\ref{ll}).

For the rest of the proof we only consider sufficiently large indices $j$ such that $x_j\not=y_j$. We now make use of Ishii's lemma as follows: Since $(x_j,y_j)$ is a maximum of $\psi_j(x,y)$, there are matrices $X_j,Y_j\in \mathcal{S}^{n\times n}(\R)$ with $X_j\le Y_j$ such that
\[(\eta_j,X_j)\in \overline{J}^{2,+}u(x_j),\;\;\;\;\;(\eta_j,Y_j)\in  \overline{J}^{2,-}v(y_j).\]
In fact,
\[\eta_j=D_x\left(\frac j4|x-y|^4\right)=-D_y\left(\frac j4|x-y|^4\right)=j|x_j-y_j|^2(x_j-y_j)\neq o.\]

\noindent Since $u$ is a subsolution and $v$ is a supersolution we have the following inequalities, where we write $\eta'_j:=|\eta_j|^{-1}\eta_j$:
\begin{align*}
	f(u(x_j))+g(u(x_j))\varpi(|\eta_j|)&\le \left\langle X_j\eta'_j,\eta'_j\right\rangle\\[.3cm]
	&\le \left\langle Y_j\eta'_j,\eta'_j\right\rangle\\[.3cm]
	&\le f(v(y_j))+g(v(y_j))\varpi(|\eta_j|).
\end{align*}

\noindent Thus,  for sufficiently large $j$, we find
\begin{equation}\label{li}
	\left(f(u(x_j))-f(v(y_j))\right)+\left(g(u(x_j))-g(v(y_j))\right)\varpi(|\eta_j|)\le 0.
\end{equation}
Since $u(x_j)>v(y_j)$ by (\ref{ll}), condition
$\bm{(\mathcal P)\,(\mathbf b)}$ shows that (\ref{li}) is impossible. This
concludes the proof of the proposition.\qed

\vspace{.2cm}

Suppose now that $f$ and $g$ satisfy condition
$\bm{(\mathcal P)\,(\mathbf a)}$. Given constants $a>0$ and $0<q<2$,
consider the following initial-value problem:
\begin{equation}\label{ivp}
	\left\{\begin{array}{lcl}
		\phi''(r)&\!\!\!=&\!\!\!f(\phi)+g(\phi)|\phi'|^q\;\;\;\text{for}\;r>0\\[.1cm]
		\phi(0)&\!\!\!=&\!\!\!a,\;\;\phi'(0)=0.\tag{\textup{IVP}$(a)$}
	\end{array}
	\right.
\end{equation}
Problem (\ref{ivp}) admits a solution $\vp\in C^2([0,R))$ for some
$R:=R(a)$. Local existence follows from the discussion preceding
\cite[Lemma~2.1]{CLV1}; strict increase and convexity follow from that lemma
with $c=1$.

\vspace{.2cm}
\noindent For future reference we also note the following.

\begin{lem}\label{op} Assume that $f$ and $g$ satisfy condition $\bm{(\mathcal P)\,(\mathbf a)}$.
	Given a constant $a>0$, let $\vp$ be a solution of (\ref{ivp}) in an interval $[0,R)$. If $w(x):=\vp(|x-z|)$ for some $z\in\R^n$, then $w$ is a viscosity solution of
	\begin{equation}\label{fb}
		\DN w=f(w)+g(w)|Dw|^q\;\;\;\text{in}\;\;B:=B(z,R).
	\end{equation}
	
\end{lem}

\noindent{\bfseries Proof:}  Since $w\in C^2(B\setminus\{z\})$, it is easily seen that $w$ is a classical solution of the PDE in the punctured ball $B\setminus\{z\}$.  So it suffices to show that $w$ is a viscosity solution of (\ref{fb}) at $x=z$. Since $\vp'(0)=0$ we see that $Dw(z)=o$. Now, suppose for some $\psi\in C^2(B)$ the function $w-\psi$ has a local maximum at $z$. Then $D\psi(z)=Dw(z)=o$. Note that $w(x)-w(z)\le \psi(x)-\psi(z)$ in a neighborhood of $z$. Therefore, as $x\to z$ we have
\begin{equation}\label{ts}
	w(x)-w(z)\le \psi(x)-\psi(z)=\frac12\langle D^2\psi(z)(x-z),x-z\rangle+o(|x-z|^2).
\end{equation}
Let $x=z+te$ for $t>0$, where $|e|=1$. Using this in (\ref{ts}) we find, as $t\to 0$,
\[\vp(t)-\vp(0)\le \frac{t^2}{2}\left\langle D^2\psi(z)e,e\right\rangle+o(t^2).\]

\noindent Dividing through by $t^2$, and then taking the limit as $t\to0^+$ we find
\[\frac12\vp''(0)=\lim_{t\to 0^+}\frac{\vp(t)-\vp(0)}{t^2}\le \frac12\left \langle D^2\psi(z)e,e\right\rangle.\]
By recalling that $\vp$ is a solution of (\ref{ivp}), we have
\[\frac12 f(a)\le \frac12\langle D^2\psi(z)e,e\rangle.\]
Thus we conclude that for any $e\in\mathbb{R}^n$ with $|e|=1$ we have
\[\langle D^2\psi(z)e,e\rangle\ge f(a)=f(w(z))+g(w(z))|D\psi(z)|^q.\]
Consequently
\[\Delta^{N,+}_{\infty}\psi(z)\ge f(w(z))+g(w(z))|D\psi(z)|^q.\] 

\vspace{.1cm}
Similarly, suppose that $w-\psi$ has a minimum at $z$ so that $w(x)-w(z)\ge \psi(x)-\psi(z)$ in a neighborhood of $z$. Then $D\psi(z)=Dw(z)=o$, and 
\[o(|x-z|^2)+\frac12\langle D^2\psi(z)(x-z),x-z\rangle =\psi(x)-\psi(z)\le w(x)-w(z).\]
Given $e\in\mathbb{R}^n$ such that $|e|=1$, we take $x=z+te$, for $t>0$. Then
we find
\[\frac{o(t^2)}{t^2}+\frac12\langle D^2\psi(z)e,e\rangle\le \frac{\vp(t)-\vp(0)}{t^2}.\]
Letting $t\to 0$ we get
\[\frac12\langle D^2\psi(z)e,e\rangle\le \frac12\vp''(0).\]
In conclusion, we have shown that
\[\langle D^2\psi(z)e,e\rangle\le \vp''(0)= f(\vp(0))+g(\vp(0))|\vp'(0)|^q.\] Hence
\[\Delta^{N,-}_{\infty}\psi(z)\le f(w(z))+g(w(z))|D\psi(z)|^q.\] \qed

Our next goal is to derive a global upper bound on all viscosity subsolutions
of (\ref{me}). For this we need suitable growth conditions on $f$ and $g$ at
infinity. This growth is captured by the {\bfseries (KO)$_q$} condition,
which will allow us to show that
\[u(x)\le \mathcal Q(d_\Omega(x))\]
for some nonincreasing function
$\mathcal Q:(0,\infty)\to[0,\infty)$. Then we can write the equation,
assuming that $u>0$, as
\[\DN u=\left(\frac{f(u)}{u}\right)u+\left(\frac{g(u)}{u^{1-q}}\right)|Du|^q u^{1-q}.\]
Define
\[A(x):=\frac{f(u(x))}{u(x)},\qquad
B(x):=\frac{g(u(x))}{u(x)^{1-q}}.\]
We will require additional conditions on $f$ and $g$ so that
$0\le A(x)\le A_0$ and $0\le B(x)\le B_0$ for some positive constants
$A_0$ and $B_0$. We then apply Theorem \ref{Hi-1}.

From this point onward, assume that $f$ and $g$ satisfy
$\bm{(\mathcal P^+)}$. With this goal in mind, we study solutions of the
initial-value problem (\ref{ivp}). Given $a>0$ and $0<q<2$, let
$\phi\in C^2([0,R))$ be a solution of (\ref{ivp}), where $0<R\le\infty$ and
$[0,R)$ is the maximal interval of existence. To emphasize its dependence on
the initial value $a>0$, we also write $R$ as $R(a)$. Since
$\phi(r)\ge a>0$, condition $\bm{(\mathcal P^+)}$ gives
\[
\phi''(r)=f(\phi(r))+g(\phi(r))\phi'(r)^q>0,
\qquad \phi'(r)>0\quad (0<r<R(a)).
\]
From the equation in (\ref{ivp}) we find that
\[\phi''\phi'\ge f(\phi)\phi',\;\;\text{and}\;\;\;\phi''(\phi')^{1-q}\ge g(\phi)\phi'.\]
Given $0<r<R$, we integrate each of these on $(0,r)$ and we obtain
\[\phi'(r)\ge \sqrt{\mathcal F(\phi(r),a)},\;\;\;\text{and}\;\;\;\phi'(r)\ge \left((2-q)\mathcal G(\phi(r),a)\right)^{\frac{1}{2-q}}.\]

\noindent Here, for $0<t<s$ we used
\[\mathcal F(s,t):=F(s)-F(t),\;\;\text{and}\;\;\;\mathcal G(s,t):=G(s)-G(t).\]

\noindent Thus, we have

\begin{equation}\label{lo}
C(q)\le 	\frac{\phi'(r)}{\sqrt{\mathcal F(\phi(r),a)}+\mathcal G(\phi(r),a)^{\frac{1}{2-q}}},
\end{equation}
where $C(q)$ is a positive constant that depends only on $q$. Integrating (\ref{lo}) on $(0,r)$ for any $0<r<R$, we find

\begin{equation}\label{ul}
	C(q)r\le \int_a^{\phi(r)}\frac{ds}{\sqrt{\mathcal F(s,a)}+\mathcal G(s,a)^{\frac{1}{2-q}}}.
\end{equation}

\vspace{.2cm}
\noindent Assume, in addition to $\bm{(\mathcal P^+)}$, that condition
{\bfseries (KO)$_q$} holds, and define $\Psi:(0,\infty)\to(0,\infty)$ by
\begin{equation}\label{invpsi}
	\Psi(t):=\frac{1}{C(q)}\int_t^\infty\frac{ds}{\sqrt{\mathcal F(s,t)}+\mathcal G(s,t)^{\frac{1}{2-q}}}.
\end{equation} 
\vspace{.2cm}

For the applications of \cite{TA-0}, define
\[
t_+:=\max\{t,0\},\qquad
\widetilde f(t):=f(t_+),\qquad
\widetilde g(t):=g(t_+).
\]
Then $\widetilde f,\widetilde g:\mathbb R\to[0,\infty)$ are continuous and
nondecreasing and agree with $f,g$ on $[0,\infty)$. The results quoted from
\cite{TA-0} are applied to $\widetilde f$ and $\widetilde g$. Since all the
integrals below have positive arguments and the solution of \textup{(IVP)}
with $\phi(0)=a>0$ remains positive, this extension does not change any of the
displayed formulas or the ODE solution under consideration.

\begin{lem}\label{psi-properties}
Suppose that $f$ and $g$ satisfy $\bm{(\mathcal P^+)}$ and {\bfseries(KO)$_q$}. Then the function $\Psi$ defined in (\ref{invpsi}) is finite, continuous, and strictly decreasing on $(0,\infty)$.
\end{lem}

\noindent{\bfseries Proof:}
Applied to the nonnegative extensions $\widetilde f$ and $\widetilde g$
defined above, \cite[Appendix B.3.1 and B.3.3,
pp.~3854--3856]{TA-0} with $k=1$ gives the finiteness and nonincreasing
monotonicity of $\Psi$.

After the change of variables $s=t+r$, set
\[
K_t(r):=\left[
\left(\int_0^r f(t+\tau)\,d\tau\right)^{1/2}
+\left(\int_0^r g(t+\tau)\,d\tau\right)^{\frac{1}{2-q}}
\right]^{-1}.
\]
Then $C(q)\Psi(t)=\int_0^\infty K_t(r)\,dr$. If $0<\tau<t$, monotonicity
of $f$ and $g$ gives $K_\tau(r)\ge K_t(r)$ for every $r>0$. Since at least
one of $f$ and $g$ is strictly increasing, this inequality is strict for
every $r>0$; hence $\Psi(\tau)>\Psi(t)$.

It remains to prove continuity. If $t_j\to t>0$, continuity
of $f$ and $g$ gives $K_{t_j}(r)\to K_t(r)$ for every $r>0$. For all large
$j$, monotonicity gives $0<K_{t_j}(r)\le K_{t/2}(r)$, and
$K_{t/2}$ is integrable because its integral equals $C(q)\Psi(t/2)$.
The dominated convergence theorem therefore proves that
$\Psi(t_j)\to\Psi(t)$.
\qed

\begin{rem}\label{ull}
Suppose that $f$ and $g$ satisfy $\bm{(\mathcal P^+)}$ and {\bfseries(KO)$_q$}.
If $0<q\le 1$, then
\begin{equation}\label{u-l}
	\lim_{t\to\infty}\Psi(t)=0.
\end{equation}
When $1<q<2$, the limit (\ref{u-l}) may fail under {\bfseries(KO)$_q$}
alone, but it holds if {\bfseries(Q-q)} is also satisfied. These assertions follow
by taking $k=1$ in \cite[Appendix B.3.3--B.3.4,
pp.~3855--3858]{TA-0}: the integral in equation~\textup{(B.10)} of that
paper is precisely $C(q)\Psi(t)$.
\end{rem}

Under the hypotheses of Remark \ref{ull}, we have $\Psi(\infty)=0$. Hence, by Lemma \ref{psi-properties},
\[
\Psi(\mathbb R^+)=(0,\Psi(0+)),
\]
and $\Psi$ has a decreasing inverse
\[
\Phi:(0,\Psi(0+))\longrightarrow(0,\infty).
\]
When $\Psi(0+)<\infty$, we extend the inverse continuously to the endpoint by defining $\Phi(\Psi(0+))=0$. Furthermore, from (\ref{ul}), we see that
\begin{equation}\label{Ra}
R(a)\le \Psi(a),\;\;\;a>0.
\end{equation}

In particular, $R(a)<\infty$. The blow-up behavior of the maximal solution
follows by taking $n=k=1$ in
\cite[Theorem~5.2(ii), pp.~3834--3837]{TA-0}. Indeed, equation~\textup{(4.1)}
of that paper then reduces precisely to the initial-value problem
\textup{(\ref{ivp})}, and the theorem gives
\[
\phi(r)\longrightarrow\infty
\qquad\text{as }r\uparrow R(a).
\]

\noindent Next, we derive a global upper estimate for subsolutions of equation (\ref{me}). We write 
\[\mathcal Q(t):=\Phi\left(\min\left\{t,\Psi(0+)\right\}\right),\;\;\;t>0,\]
where the endpoint value of $\Phi$ is understood as above.

\begin{prop}[Global $L^\infty$ Estimate]\label{ge}
	Suppose that the functions $f$ and $g$ satisfy $\bm{(\mathcal P^+)}$ and {\bfseries(KO)$_q$}. If $0<q\le 1$, then
\begin{equation}\label{ges}
u(x)\le \mathcal Q(d_\Omega(x)),\;\;\;x\in\Om,
\end{equation} 
for any viscosity subsolution $u\in \textup{USC}(\Omega)$ of equation (\ref{me}). If, on the other hand, $1<q<2$, then estimate (\ref{ges}) holds, provided that $f$ and $g$ also satisfy condition {\bfseries(Q-q)}.
	
\end{prop}

\noindent{\bfseries Proof:} Let $x\in\Omega$ and first assume that
$0<d_\Omega(x)<\Psi(0+)$. Take any $a>\Phi(d_\Omega(x))$ and consider a
solution $\phi\in C^2([0,R(a)))$ of the initial-value problem (\ref{ivp}) on
its maximal interval of existence. By the preceding cited result,
$\phi(r)\to\infty$ as $r\uparrow R(a)$. According to (\ref{Ra}), we have
\[R(a)\le \Psi(a)<d_\Omega(x).\]
Therefore $B(x,R(a))\subset\Omega$. We will show that $u(x)\le a$, which
yields $u(x)\le\Phi(d_\Omega(x))$. To this end, let
$w(y):=\phi(|x-y|)$ for $y\in B(x,R(a))$. It follows from Lemma \ref{op}
that $w$ is a viscosity solution of (\ref{me}) in $B(x,R(a))$.

Since $R(a)<d_\Omega(x)$, the closed ball $\overline{B(x,R(a))}$ is a compact subset of $\Omega$. Hence, by the upper semicontinuity of $u$, the quantity
\[
M:=\max_{\overline{B(x,R(a))}}u
\]
is finite. Since $\phi(r)\to\infty$ as $r\uparrow R(a)$, we may choose $0<\rho<R(a)$ sufficiently close to $R(a)$ that $\phi(\rho)>M$. Therefore, for every $y\in\partial B(x,\rho)$,
\[
u(y)\le M<\phi(\rho)=w(y).
\]
By the comparison principle, Proposition \ref{cp}, applied in $B(x,\rho)$, we conclude that $u\le w$ in $B(x,\rho)$. In particular,
\[u(x)\le w(x)=\phi(0)=a.\]
Since $a>\Phi(d_\Omega(x))$ is arbitrary, we conclude that $u(x)\le \Phi(d_\Omega(x))$. 
Next, let us suppose that $d_\Omega(x)\ge \Psi(0+)$ when
$\Psi(0+)<\infty$. Then $\Psi(a)<d_\Omega(x)$ for any $a>0$, so that
$B(x,R(a))\subset\Omega$ for any $a>0$. Following the same argument used
above, we see that $u(x)\le a$ for any $a>0$. This shows that
$u(x)\le 0=\Phi(\Psi(0+))$.

Consequently,
\[u(x)\le \mathcal Q(d_\Omega(x)),\;\;x\in\Omega,\] where $\mathcal Q$ is the nonincreasing function
\[\mathcal Q(t):=\Phi\left(\min\{t,\Psi(0+)\}\right),\;\;\;t>0.\]\qed

\begin{rem}
It follows from the estimate (\ref{ges}) that if $u$ is a subsolution of
(\ref{me}), and $u(x_0)>0$ at some point $x_0\in\Omega$, then
$d_\Omega(x_0)<\Psi(0+)$, and therefore
$u(x_0)\le \Phi(d_{\Omega}(x_0))$.
\end{rem}

\vspace{.2cm}
\section{Proof of the Harnack Inequalities}
We are now ready to prove the two Harnack inequalities. The first concerns
positive supersolutions of (\ref{e-1}), and the second concerns nonnegative
solutions of (\ref{me}). The proof is an adaptation of those used in
\cite{TBH, BIV}.

\vspace{.2cm}
\noindent{\bfseries Proof of Theorem \ref{Hi-1}:} Let $\alpha>1$. Given $x\in\Omega$  and $r>0$, let
\[v_x(\xi):=\left(1-\frac{|\xi-x|}{r}\right)^\alpha,\;\;\;\;\xi\in B(x,r).\]
For a positive function $w$, we consider the operator 
\[Lw(\xi):=\Delta_\infty^N w(\xi)-A_0w(\xi)-B_0|Dw(\xi)|^qw(\xi)^{1-q}.\]
On observing that
\[|Dv_x(\xi)|=\frac\alpha r\left(1-\frac{|\xi-x|}{r}\right)^{\alpha-1},\;\;\;\text{and}\;\;\;\;\Delta_\infty^Nv_x(\xi)=\frac{\alpha(\alpha-1)}{r^2}\left(1-\frac{|\xi-x|}{r}\right)^{\alpha-2},\]
set $s:=1-|\xi-x|/r$. We then compute
\begin{align*}
Lv_x(\xi)
&=\frac{s^{\alpha-2}}{r^2}
\left[\alpha(\alpha-1)-A_0r^2s^2
-B_0\alpha^qr^{2-q}s^{2-q}\right]\\[.2cm]
&\ge \frac{s^{\alpha-2}}{r^2}
\left[\alpha(\alpha-1)-|A_0|r^2-|B_0|\alpha^qr^{2-q}\right].
\end{align*}
We choose $r_0:=r_0(\alpha,q,A_0,B_0)$ such that \[\alpha(\alpha-1)-|A_0|r^2 -|B_0|\alpha^qr^{2-q}>0\;\;\;\;\forall\;0<r<r_0.\]
In fact, we can take
\[r_0:=\min\left\{\left(\frac{\alpha(\alpha-1)}{|A_0|+|B_0|\alpha^q}\right)^{\frac12}\,,\,\left(\frac{\alpha(\alpha-1)}{|A_0|+|B_0|\alpha^q}\right)^{\frac{1}{2-q}}\right\},\]
with the usual convention that the corresponding restriction is omitted if $A_0=B_0=0$.

\noindent Thus, for  $0<r<r_0$, we have
\begin{equation}\label{pos}
Lv_x(\xi)>0\;\;\;\text{in}\;\mathcal O_x:=B(x,r)\setminus\{x\}.
\end{equation}

Let $u\in \textup{LSC}(\Om)$ be a positive viscosity supersolution of
(\ref{e-1}). Fix $0<r<r_0$ such that $B(x_0,2r)\subset\Omega$. Let
$x\in B(x_0,r/3)$ be arbitrary. We consider the barrier
\begin{equation}\label{wx}
	w_x(\xi):=u(x)v_x(\xi)\;\;\text{for}\;\;\xi\in B(x,r).
\end{equation}

\noindent We observe that $w_x(x)=u(x)$, and $w_x(\xi)=0$ for $|\xi-x|=r$. We claim that $w_x\le u$ on $\overline{B(x,r)}$. To prove the claim, let us consider the following closed set
\[S:=\left\{t\in[0,1]: u(\xi)\ge tw_x(\xi)\;\forall \xi\in \overline{B(x,r)}\right\}.\]
Since $u$ is positive and lower-semicontinuous, we note that
\[\beta:=\min_{\overline{B(x,r)}} u>0.\]
Since $x\in\overline{B(x,r)}$, we have $0<\beta\le u(x)$. Moreover,
$w_x(\xi)=u(x)v_x(\xi)\le u(x)$ in $\overline{B(x,r)}$. Hence, if
$0\le t\le\beta/u(x)$, then
\[
tw_x(\xi)\le t u(x)\le\beta\le u(\xi)
\qquad\text{for every }\xi\in\overline{B(x,r)}.
\]
Thus $[0,\beta/u(x)]\subseteq S$, so that $S$ is a nonempty closed set. Now, let 
\[t_0:=\sup S.\] 
Then $\beta/u(x)\le t_0\le 1$, and $t_0\in S$. We want to show that $t_0=1$. Assume that $t_0<1$. We take sequences $t_0<t_n<1$ and $\xi_n\in \overline{B(x,r)}$ such that $t_n\to t_0$ and 
\[u(\xi_n)<t_n w_x(\xi_n),\;\;\;\forall\,n.\] We can suppose that $\xi_n\to \xi_0\in \overline {B(x,r)}$. Then, by the lower semicontinuity of $u$ we see that
\[u(\xi_0)\le t_0w_x(\xi_0).\] Let us write $w_x^0:=t_0w_x$.  Since $t_0\in S$ and $\xi_0\in \overline{B(x,r)}$ we conclude that
\[u(\xi_0)=w_x^0(\xi_0).\]
Since $u>0$ and $w_x=0$ on $\partial B(x,r)$, we have
$\xi_0\notin\partial B(x,r)$. The case $\xi_0=x$ is impossible, because
$u(x)=w_x^0(x)=t_0u(x)$ and $u(x)>0$ would imply $t_0=1$, contrary to the
assumption that $t_0<1$. Therefore
$\xi_0\in B(x,r)\setminus\{x\}$, and hence $u-w_x^0$ has a minimum at
$\xi_0\in\mathcal O_x$.

Since $u$ is a viscosity supersolution of (\ref{e-1}) and
$w_x^0\in C^2(\mathcal O_x)$, and since $u-w_x^0$ has a minimum at
$\xi_0\in\mathcal O_x$, we have
\begin{align}\Delta_\infty^N w^0_x(\xi_0)&\le A(\xi_0)u(\xi_0)+B(\xi_0)|Dw^0_x(\xi_0)|^qu(\xi_0)^{1-q}\notag\\[.2cm]
	&\le A_0u(\xi_0)+B_0|Dw^0_x(\xi_0)|^qu(\xi_0)^{1-q}.\label{sup}
\end{align}

On the other hand, recalling that $\Delta_\infty^N$ is homogeneous of degree one, we have
\begin{align*}
\Delta_\infty^N w^0_x(\xi_0)&=t_0u(x)\Delta_\infty^N v_x(\xi_0)\\[.2cm]
&>t_0u(x)\left[A_0v_x(\xi_0)+B_0|Dv_x(\xi_0)|^qv_x(\xi_0)^{1-q}\right],\;\;\text{by}\;(\ref{pos})\\[.2cm]
&=A_0w^0_x(\xi_0)+B_0|Dw^0_x(\xi_0)|^qw^0_x(\xi_0)^{1-q}\\[.2cm]
&=A_0u(\xi_0)+B_0|Dw_x^0(\xi_0)|^qu(\xi_0)^{1-q},
\qquad\text{since }w^0_x(\xi_0)=u(\xi_0).
\end{align*}
This last inequality contradicts (\ref{sup}), thus proving our claim that $t_0=1$. This shows that 
\begin{equation*}
	u(\xi)\ge w_x(\xi),\;\;\;\;\xi\in \overline{B(x,r)}.
\end{equation*}

The above observation allows us to establish the desired Harnack inequality
as follows. Let $x_0\in\Omega$, and fix $0<r<r_0$ such that
$B(x_0,2r)\subseteq\Omega$. We take an arbitrary point
$y\in B(x_0,r/3)$.

Note that
\[ y\in B(x,2r/3)\subset B(x,r)\subseteq B(x_0,2r).\] With $w_x$ defined as in (\ref{wx}), we estimate
\begin{align*}
	u(y)\ge w_x(y)&=u(x)\left(1-\frac{|y -x|}{r}\right)^\alpha\\[.3cm]
	&\ge \frac{1}{3^\alpha}u(x).
\end{align*}
This completes the proof.\qed

\vspace{.3cm}

We need one more observation to prepare for the proof of Theorem \ref{mhi}.
Let the function $h$ defined in {\bfseries(f-g)} satisfy {\bfseries(Q-m)}
and {\bfseries(KO)$_q$}. Then
\[\frac{F(s)^{\frac12}+(G(s))^{\frac{1}{2-q}}}{s}\le h(s)\le \frac{h(1)}{\mu},\;\;\;\;\;0<s\le 1.\]
This leads to the following integral inequality:
\[\frac{\mu}{h(1)}\int_t^1\frac{ds}{s}\le \int_t^\infty\frac{ds}{(F(s))^{\frac12}+(G(s))^{\frac{1}{2-q}}}\le C\Psi(t),\;\;\;\;0<t\le 1.\]
As a consequence of this we see that

\begin{equation}\label{Psi-0}
	\Psi(0+):=\lim_{t\to 0^+}\Psi(t)=\infty.
\end{equation}
This leads to the following remark.

\begin{rem}\label{ulep}
	If the hypotheses of Proposition \ref{ge} hold and the function $h$ given
	in {\bfseries(f-g)} satisfies {\bfseries(Q-m)}, then $\Psi(0+)=\infty$ and the
	inequality in (\ref{ges}) becomes
	\[u(x)\le \Phi(d_\Omega(x))\;\;\;\;\forall\,x\in\Omega.\]
\end{rem}

Now, we are ready for the proof of Theorem \ref{mhi}.

\vspace{.1cm}

\noindent{\bfseries Proof of Theorem \ref{mhi}:} Fix $\alpha:=2$, and let
$u\in\mathcal C(\Omega)$ be a nonnegative viscosity solution of (\ref{me}).
Under the assumptions given in {\bfseries(1)} or {\bfseries(2)}, it follows that the
Keller--Osserman condition {\bfseries(KO)$_q$} holds. Therefore,
Proposition \ref{ge} applies and gives
\[u(x)\le \mathcal Q(d_\Omega(x))\;\;\;\;\;\forall x\in\Omega.\] Moreover, {\bfseries(Q-M)} implies {\bfseries(Q-m)} in case {\bfseries(1)}, while {\bfseries(Q-m)} is assumed in case {\bfseries(2)}. Hence Remark \ref{ulep} applies in both cases, and the estimate improves to
\begin{equation}\label{gest}
	u(x)\le \Phi(d_\Omega(x))\;\;\;\;\forall x\in\Omega.
\end{equation}
Let
\[\delta:=\frac13\operatorname{dist}(\mathcal O,\partial\Omega),
\qquad \Omega_0:=\{x\in\Omega:d_\Omega(x)>\delta\}.\]
Then $\mathcal O\subset\Omega_0\subset\Omega$.

From (\ref{gest}) we infer that
\begin{equation}\label{do}
	0\le u\le \Phi(\delta)\;\;\;\text{in}\;\Omega_0.
\end{equation}
Let $0<\varepsilon <\Phi(\delta)$  and set
\[w_\varepsilon:=u+\varepsilon.\]
Then $w_\varepsilon>0$ and is a positive viscosity solution of 
\[\Delta_\infty^N w_\varepsilon=A_\varepsilon(x)w_\varepsilon+B_\varepsilon (x)|Dw_\varepsilon|^q|w_\varepsilon|^{1-q},\]
where
\[A_\varepsilon (x):=\frac{f(u(x))}{u(x)+\varepsilon},\;\;\text{and}\;\;\;B_\varepsilon (x):=\frac{g(u(x))}{(u(x)+\varepsilon)^{1-q}},\;\;\;x\in\Omega.\]
We begin with the proof of {\bfseries(1)} of the theorem, where $0<q<1$. In this case, the goal is to find a uniform estimate of
\[\sqrt{A_\varepsilon(x)}+B_\varepsilon(x)^{\frac{1}{2-q}}\]
by constants $A_0$ and $B_0$ that depend on $q,\nu,\varrho,\mathcal O$, and
$\operatorname{dist}(\mathcal O,\partial\Omega)$ only. In particular, $A_0$
and $B_0$ are independent of $\varepsilon$ and $u$. By (\ref{do}) and
Lemma \ref{equasimon}, applied with $s=\Phi(\delta)$ and $t=u(x)$, we have
for $x\in\Omega_0$
\begin{align}
\sqrt{A_\varepsilon(x)}+B_\varepsilon(x)^{\frac{1}{2-q}}
&=h_\varepsilon(u(x))\notag\\
&\le \frac{2}{\varrho}
\frac{\Phi(\delta)+\varepsilon}{\Phi(\delta)}
h_\varepsilon(\Phi(\delta))\notag\\
&\le \frac{4}{\varrho}h_\varepsilon(\Phi(\delta)),
\qquad 0<\varepsilon\le\Phi(\delta).\label{ho}
\end{align}
In this case $h_\varepsilon\le h$. Hence (\ref{ho}) and Lemma \ref{ap-2} give
\begin{align*}
\sqrt{A_\varepsilon(x)}+B_\varepsilon(x)^{\frac{1}{2-q}}&\le 
\frac{4}{\varrho}h_\varepsilon(\Phi(\delta))
\le\frac{4}{\varrho}h(\Phi(\delta))\\[.2cm]
&\le \frac{C(q,\nu,\varrho)}{\delta},\;\;\;\text{by Lemma \ref{ap-2}}
\end{align*}

Now, let us suppose that the assumptions in {\bfseries(2)} hold. In this case,
$1\le q<2$, and condition {\bfseries(Q-m)} states that
\[
h(s)\ge\mu h(t),\qquad s\ge t>0.
\]
We therefore apply Lemma \ref{equasimon} with its constant $\varrho$
replaced by $\mu$. Since $0\le u(x)\le\Phi(\delta)$ in $\Omega_0$, the
lemma, with $s=\Phi(\delta)$ and $t=u(x)$, gives
\begin{align*}
\sqrt{A_\varepsilon(x)}+B_\varepsilon(x)^{\frac{1}{2-q}}
&=h_\varepsilon(u(x))\\
&\le
\frac{2}{\mu}
\frac{\Phi(\delta)+\varepsilon}{\Phi(\delta)}
h_\varepsilon(\Phi(\delta))\\
&\le
\frac{4}{\mu}
\left[
\left(\frac{f(\Phi(\delta))}{\Phi(\delta)}\right)^{\frac12}
+
\left(2\Phi(\delta)\right)^{\frac{q-1}{2-q}}
g(\Phi(\delta))^{\frac{1}{2-q}}
\right]\\
&=:C(\mu,q,f,g,\delta),
\qquad x\in\Omega_0,
\end{align*}
where we used $0<\varepsilon\le\Phi(\delta)$ in the last inequality.

In summary, we have shown that there are constants $\tilde A_0$ and $\tilde B_0$ such that
\[0\le A_\varepsilon(x)\le \tilde A_0, \;\;\;\text{and}\;\;\;0\le B_\varepsilon(x)\le \tilde B_0\;\;\;\;\text{in}\;\;\Omega_0.\]
We note that, when $0<q<1$, these depend only on $q,\nu,\varrho$ and
$\operatorname{dist}(\mathcal O,\partial\Omega)$. On the other hand, these
coefficients depend on $q,\mu,f,g,\mathcal O$ and
$\operatorname{dist}(\mathcal O,\partial\Omega)$ when $1\le q<2$.

Let $r_0$ be the constant in Theorem~\ref{Hi-1}, corresponding to the
bounds $\tilde A_0$ and $\tilde B_0$, and choose
\[
0<r<\min\left\{\frac23\delta,r_0\right\}.
\]
For $K:=\overline{\mathcal O}$, this choice ensures that
$B(z,2r)\subset\Omega_0$ for every $z\in K$. By compactness and
connectedness of $K$, a standard finite-cover argument yields an integer
$m=m(K,r)$ such that any two points of $K$ can be joined by a chain of at
most $m$ overlapping balls $B(z_j,r/3)$, where $z_j\in K$ and
$B(z_j,2r)\subset\Omega_0$. For example, one may take a finite cover of $K$
by the balls $B(z,r/6)$ and use the connected overlap graph.

Applying Theorem~\ref{Hi-1} with $\alpha=2$ to
$w_\varepsilon=u+\varepsilon$ on each ball in the chain gives the local
factor $3^2=9$. Iteration across the overlaps therefore yields
\[
w_\varepsilon(x)\le 9^m w_\varepsilon(y)
\qquad\text{for all }x,y\in K.
\]
Choosing $x$ and $y$ at which $u$ attains its maximum and minimum on $K$,
respectively, and then letting $\varepsilon\to0^+$, we obtain
\[
\sup_{\mathcal O}u\le 9^m\inf_{\mathcal O}u.
\]
This proves \eqref{Hi-2}.\qed

\vspace{.5cm}

\appendix

\section{Auxiliary Growth Estimates}\label{app:growth}

In this appendix we will prove many of the technical results involving $f$ and $g$ that were used in the proof of the Harnack inequality for (\ref{me}).
We start with a general discussion of positive functions on $\mathbb R^+$ that satisfy condition {\bfseries(Q-M)}.

\begin{lem}\label{ap-1}
Let $h:(0,\infty)\to(0,\infty)$ be a continuous function that satisfies
condition {\bfseries(Q-M)}. Then

\vspace{.1cm}
\begin{enumerate}[\bfseries(a)]
	\item $\displaystyle{\lim_{t\to\infty}\frac{t^\eta}{h(t)}=0}$ for any $0<\eta<\nu$.
	
	\vspace{.2cm}
	
\item $\displaystyle{\int_t^\infty\frac{ds}{sh(s)}
\le \frac{1}{\nu\varrho h(t)}}$ for every $t>0$.
\end{enumerate}
\end{lem}
\noindent {\bfseries Proof:} 

\noindent{\bfseries(a)} Condition {\bfseries(Q-M)} shows that for any $t\ge 1$,
\[0\le \frac{t^\eta}{h(t)}\le\frac{1}{\varrho h(1)} \frac{1}{t^{\nu-\eta}}.\]
Since $0<\eta<\nu$, the desired limit follows from this.

\vspace{.2cm}
\noindent{\bfseries(b)} Given $t>0$, we see that for all $s\ge t$ we have
\[sh(s)\ge \varrho s^{1+\nu}\frac{h(t)}{t^\nu}.\]
Therefore, since $\nu>0$,
\begin{align}
\int_t^\infty\frac{ds}{sh(s)}&\le \frac{t^\nu}{\varrho h(t)}\int_t^\infty\frac{ds}{s^{\nu+1}}\notag\\[.2cm]
&=\frac{1}{\nu\varrho h(t)},\;\;\;t>0.\label{h-estimate}
\end{align}
\qed

\begin{rem}
It can easily be seen that condition {\bfseries(Q-M)} implies the existence of constants $1<\rho\le \theta$ such that
\[h(\theta s)\ge \rho h(t)\;\;\;\;\forall\,s\ge t>0.\]
\end{rem}

\noindent Now let $f,g:\mathbb R\to\mathbb R$ be continuous functions that are nonnegative on $\mathbb R^+_0$ and consider the function $h$ defined in (\ref{hfg}).

Observe that
\[sh(s)=(sf(s))^{\frac12}+(sg(s))^{\frac{1}{2-q}},\;\;s>0.\]
\begin{rem}
If $f$ and $g$ are nondecreasing functions such that $h$ satisfies
{\bfseries(Q-M)}, then the pair $(f,g)$ satisfies the Keller--Osserman
condition {\bfseries(KO)$_q$}. Indeed, monotonicity gives
\[
F(2s)\ge s f(s),\qquad G(2s)\ge s g(s),qquad s>0.
\]
Consequently, part~{\bfseries(b)} of Lemma~\ref{ap-1} yields
\begin{align*}
\int_2^\infty\frac{dt}{\sqrt{F(t)}+G(t)^{\frac{1}{2-q}}}
&=2\int_1^\infty
\frac{ds}{\sqrt{F(2s)}+G(2s)^{\frac{1}{2-q}}}\\
&\le 2\int_1^\infty\frac{ds}{s h(s)}<\infty.
\end{align*}
\end{rem}

\begin{lem}\label{ap-2}
Let $f,g:\mathbb R\to\mathbb R$ satisfy condition
$\bm{(\mathcal P^+)}$, and assume that the function $h$ given in
(\ref{hfg}) satisfies condition {\bfseries(Q-M)}.
If $0<q<1$, then there is a constant $C=C(q,\nu,\varrho)>0$ such that
\begin{equation}\label{ues}
	h(\Phi(r))\le \frac Cr,\;\;\;r>0.
\end{equation}
\end{lem}
\noindent{\bfseries Proof:} The preceding remark first shows that
{\bfseries(KO)$_q$} holds. Hence Lemma \ref{psi-properties},
Remark \ref{ull}, and (\ref{Psi-0}) show that
$\Phi:(0,\infty)\to(0,\infty)$ is well defined. Since $h$ satisfies
condition {\bfseries(Q-M)}, Lemma \ref{ap-1} applies. For $s\ge 2t$ we see that
\[\mathcal F(s,t)\ge \frac12 F(s),\qquad
\mathcal G(s,t)\ge \frac12 G(s).\]
Therefore, for $t>0$,
\begin{align}
&\int_{2t}^\infty\frac{ds}{\sqrt{\mathcal F(s,t)}+\left(\mathcal G(s,t)\right)^{\frac{1}{2-q}}}\le \int_{2t}^\infty \frac{ds}{2^{-\frac12}\sqrt{ F(s)}+2^{-\frac{1}{2-q}}\left(G(s)\right)^{\frac{1}{2-q}}}\notag\\[.3cm]
&\le 2^{\frac{3-q}{2-q}}\int_t^\infty \frac{ds}{\sqrt{ F(2s)}+\left(G(2s)\right)^{\frac{1}{2-q}}}\le4\int_t^\infty\frac{ds}{sh(s)}\notag\\
&\le \frac{4}{\nu\varrho h(t)}
\qquad\text{by part {\bfseries(b)} of Lemma \ref{ap-1}.}\label{2tt}
\end{align}

\vspace{.2cm}
\noindent Now, for $0<t\le s\le 2t$, we have the following:
\begin{align}
\sqrt{\mathcal F(s,t)}+\left(\mathcal G(s,t)\right)^{\frac{1}{2-q}}&\ge \sqrt{f(t)(s-t)}+\left(g(t)(s-t)\right)^{\frac{1}{2-q}}\notag\\[.3cm]
&=\sqrt{tf(t)\left(\frac st-1\right)}+\left[tg(t)\left(\frac st-1\right)\right]^{\frac{1}{2-q}}\notag\\
&\ge th(t)\left(\frac st-1\right)^{\frac{1}{2-q}},\;\;\text{since $q\ge 0$}.\label{sq}
\end{align}

\noindent If $0<q<1$, then (\ref{sq}) shows that
\begin{equation}\label{t}
	\int_t^{2t}\frac{ds}{\sqrt{\mathcal F(s,t)}+\left(\mathcal G(s,t)\right)^{\frac{1}{2-q}}}\le \left(\frac{2-q}{1-q}\right)\frac{1}{h(t)},\;\;\;\;t>0.
\end{equation}

\noindent From (\ref{2tt}) and (\ref{t}), we conclude that the following holds for all $t>0$:
\begin{align}
\Psi(t)&=\frac{1}{C(q)}\int_t^\infty\frac{ds}{\sqrt{\mathcal F(s,t)}+\left(\mathcal G(s,t)\right)^{\frac{1}{2-q}}}\notag\\[.3cm]
&\le \frac{C}{h(t)},\qquad 0<q<1.
 \label{ps}
\end{align}
Here $C$ is a positive constant that depends only on $q,\nu$, and $\varrho$.
\vspace{.2cm}
\noindent Then, for $r>0$, estimate (\ref{ps}) with $t=\Phi(r)$ shows that
\[r\le C\frac{\Phi(r)}{(\Phi(r)f(\Phi(r)))^{\frac12}
+(\Phi(r)g(\Phi(r)))^{\frac{1}{2-q}}}.\]
In other words, for $r>0$,
\begin{equation}\label{bin}
	\frac{r\left[(\Phi(r)f(\Phi(r)))^{\frac12}+(\Phi(r)g(\Phi(r)))^{\frac{1}{2-q}}\right]}{\Phi(r)}\le C.
\end{equation}
We rewrite (\ref{bin}) as
\begin{equation}\label{uess}
h(\Phi(r))
=\left(\frac{f(\Phi(r))}{\Phi(r)}\right)^{\frac12}
+\left(\frac{g(\Phi(r))}{\Phi(r)^{1-q}}\right)^{\frac{1}{2-q}}
\le \frac Cr,\qquad r>0.
\end{equation}
Here $C$ is a positive constant that depends only on $q,\nu$, and
$\varrho$. This completes the proof of Lemma \ref{ap-2}.\qed

\begin{lem}\label{equasimon}
	Let $f,g:\mathbb R\to \mathbb R$ be nondecreasing and nonnegative on $[0,\infty)$. Let $0<q<2$, and define
	\[
	h(s):=\left(\frac{f(s)}{s}\right)^{1/2}
	+\left(\frac{g(s)}{s^{1-q}}\right)^{\frac1{2-q}},
	\qquad s>0.
	\]
	Assume that there exists $0<\varrho\le 1$ such that
	\[
	h(s)\ge \varrho h(t)
	\qquad \text{for all } s\ge t>0.
	\]
	For $\varepsilon>0$, define, for $s\ge0$,
	\[
	h_\varepsilon(s)
	:=
	\left(\frac{f(s)}{s+\varepsilon}\right)^{1/2}
	+
	\left(
	\frac{g(s)}{(s+\varepsilon)^{1-q}}
	\right)^{\frac1{2-q}}.
	\]
	Then
	\[
	h_\varepsilon(s)
	\ge
	\frac{\varrho}{2}
	\left(\frac{s}{s+\varepsilon}\right)
	h_\varepsilon(t)
	\qquad \text{for all } s\ge t\ge0.
	\]
\end{lem}

\noindent{\bfseries Proof:}
Set
\[
\beta:=\frac{1-q}{2-q},
\qquad
x_\tau:=\frac{\tau}{\tau+\varepsilon},
\qquad \tau>0.
\]
Then
\begin{align*}
h_\varepsilon(\tau)
&=x_\tau^{1/2}\left(\frac{f(\tau)}{\tau}\right)^{1/2}
+x_\tau^\beta
\left(\frac{g(\tau)}{\tau^{1-q}}\right)^{\frac1{2-q}}.
\end{align*}
Since $0<x_\tau<1$, $1/2\le1$, and $\beta\le1$, we have
\[
h_\varepsilon(\tau)\ge x_\tau h(\tau).
\]
Consequently, for $s\ge t>0$, the hypothesis gives
\begin{equation}\label{heps-basic}
h_\varepsilon(s)
\ge x_s h(s)
\ge \varrho x_s h(t).
\end{equation}

\noindent Suppose first that $0<q\le1$. Then $\beta\ge0$, and therefore
\[
h_\varepsilon(t)\le h(t).
\]
It follows immediately from (\ref{heps-basic}) that
\[
h_\varepsilon(s)
\ge \varrho x_s h_\varepsilon(t)
\ge \frac{\varrho}{2}x_s h_\varepsilon(t).
\]

\noindent Now suppose that $1<q<2$, so that $\beta<0$. Write
\[
h_\varepsilon(t)=F_\varepsilon(t)+G_\varepsilon(t),
\]
where
\[
F_\varepsilon(t):=\frac{f(t)^{1/2}}{(t+\varepsilon)^{1/2}},
\qquad
G_\varepsilon(t):=\frac{g(t)^{1/(2-q)}}{(t+\varepsilon)^\beta}.
\]
Since $F_\varepsilon(t)\le h(t)$, (\ref{heps-basic}) yields
\[
h_\varepsilon(s)\ge \varrho x_sF_\varepsilon(t).
\]
On the other hand, the function $G_\varepsilon$ is nondecreasing: indeed, $g$ is nondecreasing and, since $\beta<0$, the function
$\tau\mapsto(\tau+\varepsilon)^{-\beta}$ is nondecreasing. Hence
\[
h_\varepsilon(s)\ge G_\varepsilon(s)\ge G_\varepsilon(t)\ge x_sG_\varepsilon(t).
\]
Using
\[
\max\{\varrho a,b\}\ge\frac{\varrho}{2}(a+b),
\qquad a,b\ge0,
\]
we conclude that
\[
h_\varepsilon(s)
\ge\frac{\varrho}{2}x_s
\bigl(F_\varepsilon(t)+G_\varepsilon(t)\bigr)
=\frac{\varrho}{2}
\left(\frac{s}{s+\varepsilon}\right)h_\varepsilon(t).
\]

\noindent It remains to consider $t=0$. If $s=0$, the desired inequality
is immediate because $s/(s+\varepsilon)=0$. We therefore assume that
$s>0$. We first claim that $f(0)=0$. Indeed, if $f(0)>0$, then the
monotonicity of $f$ gives
\[
h(t)\ge \left(\frac{f(0)}{t}\right)^{1/2}
\qquad\text{for }0<t\le s.
\]
Letting $t\to0^+$ contradicts the hypothesis
$h(s)\ge\varrho h(t)$. Thus $f(0)=0$.

\noindent Suppose first that $0<q<1$, so that $\beta>0$. The same argument
shows that $g(0)=0$: otherwise,
\[
h(t)\ge \frac{g(0)^{1/(2-q)}}{t^\beta}
\longrightarrow\infty
\qquad\text{as }t\to0^+,
\]
again contradicting $h(s)\ge\varrho h(t)$. Consequently,
$h_\varepsilon(0)=0$, and the desired inequality follows.

\noindent If $q=1$, then $\beta=0$ and
$h_\varepsilon(0)=g(0)$. Since $g$ is nondecreasing,
\[
h_\varepsilon(s)\ge g(s)\ge g(0)=h_\varepsilon(0)
\ge \frac{\varrho}{2}\frac{s}{s+\varepsilon}h_\varepsilon(0).
\]

\noindent Finally, suppose that $1<q<2$, so that $\beta<0$. Since
$f(0)=0$, we have
\[
h_\varepsilon(0)=\frac{g(0)^{1/(2-q)}}{\varepsilon^\beta}.
\]
Moreover, the monotonicity of $g$ and the fact that $-\beta>0$ give
\begin{align*}
h_\varepsilon(s)
&\ge
\frac{g(s)^{1/(2-q)}}{(s+\varepsilon)^\beta}\\
&\ge
\frac{g(0)^{1/(2-q)}}{\varepsilon^\beta}
=h_\varepsilon(0)
\ge
\frac{\varrho}{2}\frac{s}{s+\varepsilon}h_\varepsilon(0).
\end{align*}
This proves the conclusion also when $t=0$.

This completes the proof.
\qed

\section{Sharpness of the Exponent Range}\label{app:sharpness}

We record the counterexample establishing that the range $0<q<2$ in
Theorem~\ref{Hi-1} is sharp.

\begin{rem}\label{sharp-q}
	The conclusion of Theorem~\ref{Hi-1} fails in general when $q \ge 2$.
	To see this, let $\Omega=\mathbb{R}^n$ and consider
	\begin{equation}\label{coe}
	\Delta_\infty^N u = |Du|^q u^{1-q}
	\qquad \text{in } \mathbb{R}^n,
	\end{equation}
	which corresponds to \eqref{e-1} with $A \equiv 0$ and $B \equiv 1$.

	For a positive integer $k$, define
	\[
	u_k(x)=e^{kx_1},
	\qquad x=(x_1,\dots,x_n)\in \mathbb{R}^n.
	\]
	Then $u_k \in C^\infty(\mathbb{R}^n)$ and
	\[
	Du_k = k e^{k x_1} e_1,
	\qquad
	D^2u_k = k^2 e^{k x_1}\, e_1\otimes e_1.
	\]
	Consequently,
	\[
	\Delta_\infty^N u_k
	=\frac{\langle D^2u_k Du_k, Du_k\rangle}{|Du_k|^2}
	=k^2 e^{k x_1}.
	\]
	On the other hand,
	\[
	|Du_k|^q u_k^{1-q}
	=(k e^{k x_1})^q(e^{k x_1})^{1-q}
	=k^q e^{k x_1}.
	\]
	Therefore,
	\[
	\Delta_\infty^N u_k
	=k^2 e^{kx_1}
	\le k^q e^{kx_1}
	=|Du_k|^q u_k^{1-q}
	\]
	whenever $q\ge2$. Thus $u_k$ is a positive classical supersolution of
	\eqref{coe}.

	Now let $r>0$. Since $u_k$ depends only on $x_1$,
	\[
	\sup_{B(0,r/3)}u_k=e^{kr/3},
	\qquad
	\inf_{B(0,r/3)}u_k=e^{-kr/3}.
	\]
	Hence, writing $B=B(0,r/3)$,
	\[
	\frac{\sup_B u_k}{\inf_B u_k}=e^{2kr/3}.
	\]
	Since this ratio tends to infinity as $k\to\infty$, no universal
	constant $C>0$ can satisfy
	\[
	\sup_{B(0,r/3)}u_k\le C\inf_{B(0,r/3)}u_k
	\]
	for all positive supersolutions. This proves that the conclusion of
	Theorem~\ref{Hi-1} cannot hold for $q\ge2$.
\end{rem}

\end{document}